\documentclass{amsart}
\usepackage{enumerate}
\usepackage{tabu}
\usepackage{longtable}
\usepackage{multirow}
\usepackage{amssymb} 
\usepackage{amsmath} 
\usepackage{amscd}
\usepackage{amsbsy}
\usepackage{comment}
\usepackage[matrix,arrow]{xy}
\usepackage{hyperref}

\DeclareMathOperator{\Rad}{Rad}

\DeclareMathOperator{\Norm}{Norm}

\DeclareMathOperator{\GL}{GL}

\DeclareMathOperator{\Gal}{Gal}
\DeclareMathOperator{\norm}{Norm}

\DeclareMathOperator{\ord}{ord}

\newcommand{\Q}{{\mathbb Q}}
\newcommand{\Z}{{\mathbb Z}}
\newcommand{\C}{{\mathbb C}}
\newcommand{\F}{{\mathbb F}}
\newcommand{\cA}{\mathcal{A}}

\newcommand{\cQ}{\mathcal{Q}}

\newcommand{\cE}{\mathcal{E}}

\newcommand{\OO}{{\mathcal O}}

\newcommand{\sS}{\mathfrak{S}}

\newcommand{\fq}{\mathfrak{q}}

\newcommand{\mP}{\mathcal{P}}

\def\mod#1{{\ifmmode\text{\rm\ (mod~$#1$)}
\else\discretionary{}{}{\hbox{ }}\rm(mod~$#1$)\fi}}

\begin {document}

\newtheorem{thm}{Theorem}
\newtheorem{lem}{Lemma}[section]
\newtheorem{prop}[lem]{Proposition}

\theoremstyle{definition}

\theoremstyle{remark}

\title[]{Perfect powers that are sums of consecutive cubes}

\author[Michael Bennett]{Michael A. Bennett}
\address{Department of Mathematics, University of British Columbia, Vancouver, B.C., V6T 1Z2 Canada}
\email{bennett@math.ubc.ca}

\author{Vandita Patel}
\address{Mathematics Institute, University of Warwick, Coventry CV4 7AL, United Kingdom}
\email{vandita.patel@warwick.ac.uk}

\author{Samir Siksek}
\address{Mathematics Institute, University of Warwick, Coventry CV4 7AL, United Kingdom}
\email{S.Siksek@warwick.ac.uk}
\thanks{
The first-named author is supported
by NSERC. 
The second-named author is supported by an EPSRC studentship.
The third-named author
is supported by 
the EPSRC {\em LMF: L-Functions and Modular Forms} Programme Grant
EP/K034383/1.
}

\date{\today}


\begin {abstract}
Euler noted the relation $6^3=3^3+4^3+5^3$
and asked for other instances of cubes that are sums
of consecutive cubes. Similar problems have been 
studied by Cunningham, Catalan, Gennochi, Lucas,
Pagliani, Cassels, Uchiyama, Stroeker and Zhongfeng Zhang.
In particular Stroeker determined all squares that
can be written as a sum of at most $50$ consecutive cubes.
We generalize Stroeker's work by determining all perfect
powers that are sums of at most $50$ consecutive cubes.
Our methods include descent, linear forms in two logarithms,
and Frey-Hellegouarch curves.
\end {abstract}
\maketitle

\section{Introduction} \label{intro}

Euler \cite[art.\ 249]{Euler}, 
in his 1770 \emph{Vollst\"{a}ndige Anleitung zur Algebra},
notes the relation
\begin{equation}\label{eqn:Euler}
6^3=
3^3+4^3+5^3
\end{equation} 
and asks for other instances of cubes that are sums of three
consecutive cubes. 
Dickson's 
\emph{History of the Theory of Numbers}
gives an extensive survey of early work on the problem
of cubes that are sums of consecutive cubes 
\cite[pp. 582--585]{Dickson},  and also
squares that are sums of consecutive cubes
\cite[pp. 585--588]{Dickson} with contributions
by illustrious names such as Cunningham,
Catalan, Gennochi and Lucas. 
Both problems
possess some parametric families of solutions;
one such family  was constructed
by Pagliani \cite{Pagliani}
in 1829 : 
\[
\left(\frac{v^5+v^3-2v}{6}\right)^3=\sum_{i=1}^{v^3} 
\left(\frac{v^4-3v^3-2v^2-2}{6} \, + \, i\right)^3,
\]
where the congruence restriction $v \equiv 2$ or $4 \mod{6}$
ensures integrality of the cubes.
Pagliani uses this to answer a challenge, posed presumably by the editor Gergonne,
of giving $1000$ consecutive cubes whose sum is a cube. Of course, the problem
of squares that are sums of consecutive cubes possesses the well-known
parametric family of solutions
\[
\left(\frac{d(d+1)}{2}\right)^2=
\sum_{i=1}^d i^3=\sum_{i=0}^d i^3.
\]
These questions have continued to be of intermittent  interest throughout a period
of over 200 years. 
For example, Lucas \cite[page 92]{Lucas} states incorrectly that
the only square expressible as a sum of three consecutive positive cubes is
\begin{equation}\label{eqn:Lucas}
6^2=1^3+2^3+3^3.
\end{equation}
Both Cassels \cite{Cassels}
and Uchiyama \cite{Uchiyama} determine
the squares that can be written as sums of three consecutive cubes
(without reference to Lucas) showing that the only solutions
in addition to \eqref{eqn:Lucas} are
\begin{equation}\label{eqn:Cassels}
0=(-1)^3+0^3+1^3, \qquad
3^2=0^3+1^3+2^3, \qquad
204^2=23^3+24^3+25^3.
\end{equation}
Lucas also states that the only square 
that is the sum of two consecutive positive cubes is $3^2=1^3+2^3$
and the only squares that are sums of $5$ consecutive non-negative cubes are 
\begin{gather*}
10^2=0^3+1^3+2^3+3^3+4^3, \qquad
15^2=1^3+2^3+3^3+4^3+5^3,\\
315^2=25^3+26^3+27^3+28^3+29^3,\qquad
2170^2=96^3+97^3+98^3+99^3+100^3,\\
2940^2=118^3+119^3+120^3+121^3+122^3.
\end{gather*}
These two claims turn out to be correct as shown by Stroeker \cite{Stroeker}.
In modern language,
the problem of which squares are expressible as the sum of $d$
consecutive cubes, reduces for any given $d \ge 2$, to the determination 
of integral points on a genus $1$ curve. Stroeker \cite{Stroeker}, 
using a (by now) standard 
method based on linear forms in elliptic logarithms, solves
this problem for $2 \le d \le 50$.

\medskip

The problem of expressing arbitrary perfect powers as a sum of $d$
consecutive cubes with $d$ small has received somewhat less attention, likely due to the fact that techniques
 for resolving such questions are of a much more recent vintage.
Zhongfeng Zhang \cite{Zhang}
showed that the only perfect powers that are sums
of three consecutive cubes are precisely those
already noted by Euler \eqref{eqn:Euler}, Lucas \eqref{eqn:Lucas}
and Cassels \eqref{eqn:Cassels}.
Zhang's approach is write the problem as 
\begin{equation} \label{zhangy}
y^n=(x-1)^3+x^3+(x+1)^3=3x(x^2+2),
\end{equation}
and apply a descent argument that reduces this to certain ternary
equations that have already been solved in the literature.

In this paper, we extend Stroeker's aforementioned work, determining
all perfect powers that are sums of $d$ consecutive cubes, with $2 \le d \le 50$. This upper bound is somewhat arbitrary as our techniques extend to essentially any fixed values of $d$. 
\begin{thm}\label{thm:main}
Let $2 \le d \le 50$. Let $\ell$ be a prime. The integral solutions to the equation
\begin{equation}\label{eqn:premain}
(x+1)^3+(x+2)^3+\cdots+(x+d)^3=y^\ell
\end{equation}
with $x \ge 1$ are given in Table~\ref{table:solutions}.
\end{thm}

\begin{table}
\centering
{
$\begin{tabu}{|c|l|}
\hline
d & (x,y,\ell)\\
\hline\hline

2 &
	\\
3 &
        ( 22, \pm 204, 2 ),
	(2,6,3)
	\\
4 &
	(10,20,3)
	\\
5 &
        ( 24, \pm 315, 2 ),
        ( 95, \pm 2170, 2 ),
        ( 117, \pm 2940, 2 )
    	\\
6 &
   	\\ 
7 &
        ( 332, \pm 16296, 2 )
   	\\ 
8 &
        ( 27, \pm 504, 2 )
	\\
9 &
        ( 715, \pm 57960, 2 )
	\\
10 &
	\\
11 &
        ( 1314, \pm 159060, 2 )
	\\
12 &
        ( 13, \pm 312, 2 )
	\\
13 &
        ( 143, \pm 6630, 2 ),
        ( 2177, \pm 368004, 2 )
	\\
14 &
	\\
15 &
        ( 24, \pm 720, 2 ),
        ( 3352, \pm 754320, 2 ),
        ( 57959, \pm 54052635, 2 )
   	\\ 
16 &
   	\\ 
17 &
        ( 8, \pm 323, 2 ),
        ( 119, \pm 5984, 2 ),
        ( 4887, \pm 1412496, 2 )
	\\
18 &
        ( 152, \pm 8721, 2 ),
        ( 679, \pm 76653, 2 )
   	\\ 
19 &
        ( 6830, \pm 2465820, 2 )
	\\
20 &
	(2,40,3), (14,70,3)
	\\
21 &
        ( 13, \pm 588, 2 ),
        ( 143, \pm 8778, 2 ),
        ( 9229, \pm 4070220, 2 )
	\\
22 &
   	\\ 
23 &
        ( 12132, \pm 6418104, 2 )
   	\\ 
24 &
	\\
25 &
        ( 15587, \pm 9742200, 2 ),
	(5,60,3)
   	\\ 
26 &
   	\\ 
27 &
        ( 19642, \pm 14319396, 2 )
   	\\ 
28 &
        ( 80, \pm 4914, 2 )
   	\\ 
29 &
        ( 24345, \pm 20474580, 2 )
	\\
30 &
   	\\ 
31 &
        ( 29744, \pm 28584480, 2 )
   	\\ 
32 &
        ( 68, \pm 4472, 2 ),
        ( 132, \pm 10296, 2 ),
        ( 495, \pm 65472, 2 )
   	\\ 
33 &
        ( 32, \pm 2079, 2 ),
        ( 35887, \pm 39081504, 2 )
   	\\ 
34 &
   	\\ 
35 &
        ( 224, \pm 22330, 2 ),
        ( 42822, \pm 52457580, 2 )
   	\\ 
36 &
   	\\ 
37 &
        ( 50597, \pm 69267996, 2 )
   	\\ 
38 &
   	\\ 
39 &
        ( 110, \pm 9360, 2 ),
        ( 59260, \pm 90135240, 2 )
   	\\ 
40 &
        ( 3275, \pm 1196520, 2 )
   	\\ 
41 &
        ( 68859, \pm 115752840, 2 )
   	\\ 
42 &
        ( 63, \pm 5187, 2 )
   	\\ 
43 &
        ( 79442, \pm 146889204, 2 )
	\\
44 &
   	\\ 
45 &
        ( 175, \pm 18810, 2 ),
        ( 91057, \pm 184391460, 2 )
	\\
46 &
	\\
47 &
        ( 103752, \pm 229189296, 2 )
	\\
48 &
        ( 63, \pm 5880, 2 ),
        ( 409, \pm 62628, 2 ),
        ( 19880, \pm 19455744, 2 ),
        ( 60039, \pm 101985072, 2 )
	\\
49 &
        ( 117575, \pm 282298800, 2 ),
	(290,1155,3)
	\\
50 &
        ( 1224, \pm 312375, 2 )
	\\
\hline
\end{tabu}$
}
\caption{The solutions to equation \eqref{eqn:premain}
with $2 \le d \le 50$, $\ell$ prime and $x \ge 1$.}
\label{table:solutions}
\end{table}
The restriction $x \ge 1$ imposed in the statement of Theorem~\ref{thm:main} is merely
to exclude a multitude of artificial solutions. Solutions with $x \le 0$ can
in fact be deduced easily, as we now explain :
\begin{enumerate}
\item[(i)] The value $x=0$ gives the ``trivial'' solutions $(x,y,\ell)=(0,d(d+1)/2,2)$, and
no solutions for odd $\ell$. Likewise the value $x=-1$ yields the trivial solutions 
$(x,y,\ell)=(-1,(d-1)d/2,2)$ and no solutions for odd $\ell$.
\item[(ii)] For odd exponents $\ell$, there is a symmetry between the solutions to \eqref{eqn:premain} :
\[
(x,y,\ell) \longleftrightarrow (-x-d-1,-y,\ell).
\]
This allows us to deduce, from Table~\ref{table:solutions} and (i), all solutions with 
$x \le -d-1$.
\item[(iii)] The solutions with $-d \le x \le -2$ lead to non-negative solutions with smaller values
of $d$ through cancellation (and possibly applying the symmetry in (ii)).
\end{enumerate}
Of course arbitrary perfect powers that are sums of at most $50$ consecutive cubes
can be deduced from our list of $\ell$-th powers with $\ell$ prime.

\bigskip

A sum of $d$ consecutive cubes can be written as
\[
(x+1)^3+(x+2)^3+\cdots+(x+d)^3
=\left(d x+ \frac{d(d+1)}{2}\right) \left(x^2+(d+1) x+\frac{d(d+1)}{2}  \right).
\]
Thus, to prove Theorem~\ref{thm:main}, we need to solve
the Diophantine equation
\begin{equation}\label{eqn:conscubes}
\left(d x+ \frac{d(d+1)}{2}\right) \left(x^2+(d+1) x+\frac{d(d+1)}{2}  \right)=y^\ell, 
\end{equation}
with $\ell$ prime and $2 \le d \le 50$. 
We find it convenient to rewrite \eqref{eqn:conscubes} as 
\begin{equation}\label{eqn:main}
d (2 x+ d+1) \left(x^2+(d+1) x+\frac{d(d+1)}{2}  \right)=2y^\ell. 
\end{equation}
We will use a descent argument together with the identity
\begin{equation}\label{eqn:identity}
4 \left( x^2+(d+1)x+\frac{d(d+1)}{2}\right)-(2x+d+1)^2=d^2-1.
\end{equation}
to reduce \eqref{eqn:main} to a family of ternary equations. The main purpose of this paper is to highlight the degree to which 
such ternary equations can, through a combination of techniques including descent, lower bounds for linear forms in logarithms, and appeal to the modularity of Galois representations, be nowadays completely and explicitly solved.

\medskip

We are grateful to the referee for careful reading of the paper and for 
suggesting several improvements.


\section{Proof of Theorem~\ref{thm:main} for $\ell=2$}\label{sec:l2}

Although Theorem~\ref{thm:main}  with $\ell=2$ follows from
 Stroeker's paper \cite{Stroeker},
we explain briefly how this can now be done with the
help of an appropriate
computer algebra package.

Let $(x,y)$ be an integral solution to \eqref{eqn:conscubes} with $\ell=2$.
Write
$X=dx$, and $Y=dy$. 
Then $(X,Y)$ is an integral point on the elliptic curve
\[
E_d \; : \; Y^2=\left(X+ \frac{d^2+d}{2}\right)
\left(X^2+ (d^2+d)X+ \frac{d^4+d^3}{2} \right).
\]
Using the computer algebra package \texttt{Magma} \cite{magma},
we determined the integral points on $E_d$ for $2 \le d \le 50$. 
For this computation,
\texttt{Magma}
applies the standard linear forms in elliptic logarithms method
\cite[Chapter XIII]{Smart}, which is the same method used by 
Stroeker (though the implementation is independent). From this we immediately recover the original
solutions $(x,y)$ to \eqref{eqn:conscubes}
with $\ell=2$, and the latter are found in our Table~\ref{table:solutions}.
We have checked that our solutions with $\ell=2$ are precisely those
given by Stroeker.

We shall henceforth restrict ourselves to $\ell \ge 3$.

\section{Proof of Theorem~\ref{thm:main} for $d=2$}\label{sec:d2}

Our method for general $d$ explained in later sections
fails for $d=2$. This is because of the presence of solutions $(x,y)=(-2,-1)$
and $(x,y)=(-1,1)$ to \eqref{eqn:premain} for all $\ell \ge 3$. In this section 
we treat the case $d=2$ separately, reducing to Diophantine equations that
have already been solved by Nagell.

\medskip

We consider the equation  \eqref{eqn:premain} with $d=2$. For convenience, let $z=x+1$.
The equation becomes 
$z^3+(z+1)^3=y^\ell$ which can be rewritten
as
\begin{equation}\label{eqn:twocubes}
(2z+1)(z^2+z+1)=y^\ell.
\end{equation}
Here $y$ and $z$ are integers and $\ell \ge 3$ is prime. 
Suppose first that $\ell=3$. This equation here  defines
a genus $1$ curve. We checked using \texttt{Magma}
that it is isomorphic to the elliptic curve $Y^2-9Y=X^3-27$
with Cremona label \texttt{27A1}, and that it has Mordell--Weil
group (over $\Q$) $\cong \Z/3\Z$. It follows that the only rational
points on \eqref{eqn:twocubes} with $\ell=3$ are the 
three obvious ones : $(z,y)=(-1/2,0)$, $(0,1)$ and $(-1,-1)$. 
These yield the solutions
$(x,y)=(-1,1)$ and  $(x,y)=(-2,-1)$ to \eqref{eqn:premain}.

We may thus suppose that $\ell \ge 5$ is prime.
The resultant of the two factors on the left-hand side of \eqref{eqn:twocubes}
is $3$ and, moreover, $9 \nmid (z^2+z+1)$. It follows that
either
$$
2z+1=y_1^\ell, \qquad z^2+z+1=y_2^\ell, \qquad y=y_1 y_2
$$
or
$$
2z+1=3^{\ell-1} y_1^n, \qquad z^2+z+1=3 y_2^\ell, \qquad y=3 y_1 y_2.
$$
Nagell \cite{Nagell} showed that the only integer solutions to
the equation $X^2+X+1=Y^n$
with $n \ne 3^k$ are the trivial ones with $X=-1$ or $0$. 
Nagell \cite{Nagell} also solved
the equation $X^2+X+1=3Y^n$ for $n>2$ showing that the only solutions
are again the trivial ones with $X=1$. Working back, we see
that the only solutions to \eqref{eqn:twocubes} with $\ell \ge 5$
are $(z,y)=(0,1)$ and $(-1,-1)$. These again give the solutions
$(x,y)=(-1,1)$ and $(-2,-1)$ to \eqref{eqn:premain}.

\section{Descent for $\ell \ge 5$}\label{sec:descent}

Let $d \ge 3$.
We consider equation \eqref{eqn:main} with exponent $\ell \ge 5$. The argument in this section
will need modification for $\ell=3$ which we carry out in Section~\ref{sec:l3}.
For a prime
$q$ we let
\begin{equation}\label{eqn:munu}
\mu_q=\ord_q(d^2-1) \; \; \mbox{ and } \; \; 
\nu_q=\ord_q(d),
\end{equation}
i.e. the largest power of $q$ dividing $d^2-1$ and $d$, respectively.
We associate to $q$ a finite subset $T_q \subset \Z^2$ as follows. 
\begin{itemize}
\item If $q \nmid d(d^2-1)$ then let $T_q=\{(0,0)\}$.
\item For $q=2$ we define
\[
T_2=\begin{cases}
\left\{(0,1-\nu_2)\right\} & \text{if $2 \mid d$} \\
\left\{(1,0), \; (\mu_2/2,1-\mu_2/2), \; (3-\mu_2,\mu_2-2) \right\} & \text{if $2 \nmid d$ and $2 \mid \mu_2$} \\
\left\{(1,0), \; (3-\mu_2,\mu_2-2) \right\} & \text{if $2 \nmid d$ and $2 \nmid \mu_2$}. 
\end{cases}
\]
\item For odd $q \mid d$, let
\[
T_q=\{(-\nu_q,0), \; (0,-\nu_q)\}.
\]
\item For odd $q \mid (d^2-1)$, let 
\[
T_q=\begin{cases}
\{(0,0), \; (-\mu_q,\mu_q), \; (\mu_q/2,-\mu_q/2) \} & \text{if $2 \mid \mu_q$},\\
\{(0,0), \; (-\mu_q,\mu_q)\} & \text{if $2 \nmid \mu_q$}.
\end{cases}
\]
\end{itemize}

We take $\cA_d$ to be the set of pairs of positive rationals $(\alpha,\beta)$
such that
$$
(\ord_q(\alpha),\ord_q(\beta)) \in T_q
$$
for all primes $q$.
It is clear that $\cA_d$ is a finite set, which is, in practice, easy to write
down for any value of $d$.

\begin{lem}\label{lem:descent}
Let $(x,y)$ be a solution to \eqref{eqn:main}
where $\ell \ge 5$ a prime.
Then there are rationals $y_1$, $y_2$ and a pair 
$(\alpha,\beta) \in \cA_d$ such that
\begin{equation}\label{eqn:descent}
2x+d+1=\alpha y_1^\ell, \qquad 
x^2+(d+1) x+ \frac{d(d+1)}{2}=\beta y_2^\ell.
\end{equation}
Moreover, if $3 \le d \le 50$ then $y_1$ and $y_2$ are integers.
\end{lem}

\noindent \textbf{Remark.} The reader will observe that the definition
of $\cA_d$ is independent of $\ell$. Thus, given $d$, the lemma provides us with a way
of carrying out the descent uniformly for all $\ell \ge 5$.

\begin{proof}
Let us first assume the first part of the lemma and deduce the second.
Using a short \texttt{Magma} script, we wrote down all possible
pairs $(\alpha,\beta) \in \cA_d$ for $3 \le d \le 50$ and checked
that 
\[
\max \{ \ord_q(\alpha), \;  \ord_q(\beta)  \} \le 4
\]
for all primes $q$.
As $x$ is an integer, we know from \eqref{eqn:descent} that
\[
\ord_q(\alpha)+ \ell \ord_q(y_1) \ge 0 \; \; \mbox{ and } \; \;  \ord_q(\beta)+ \ell \ord_q(y_2) \ge 0,
\]
for all primes $q$. Since $\ell \ge 5$, it is clear that $\ord_q(y_1) \ge 0$
and $\ord_q(y_2) \ge 0$ for all primes $q$. This proves the second part
of the lemma.


\medskip

We now prove the first part of the lemma. 
For $2x+d+1=0$ (which can only arise for odd values of $d$) we can take $y_1=0$, $y_2=1$,
\begin{equation}\label{eqn:trivsolution}
\alpha=\frac{8}{d(d^2-1)} \; \; \mbox{ and } \; \;  \beta=\frac{d^2-1}{4};
\end{equation}
it is easy to check that this particular pair $(\alpha,\beta)$
belongs to $\cA_d$. We shall henceforth suppose that $2x+d+1 \ne 0$.

\medskip

\noindent \textbf{Claim:}
Let $q$ be a prime and define
\[
\epsilon=\ord_q(2x+d+1) \; \; \mbox{ and } \; \; 
\delta=\ord_q \left(x^2+(d+1)x+ \frac{d(d+1)}{2}\right).
\]
Then
$(\epsilon,\delta) \equiv (\epsilon^\prime,\delta^\prime) \mod{\ell}$
for some $(\epsilon^\prime,\delta^\prime) \in T_q$. 

To complete the proof of Lemma \ref{lem:descent}, it is clearly enough to prove this claim. From \eqref{eqn:main}
and \eqref{eqn:identity}, the claim is certainly true if $q \nmid d(d^2-1)$,
so we may suppose that $q \mid d(d^2-1)$.
Observe that for any $q$, from \eqref{eqn:main},
\begin{equation}\label{eqn:val1}
\nu_q+\epsilon+\delta \equiv \ord_q(2) \mod{\ell} \, .
\end{equation}
Moreover, from \eqref{eqn:identity},
\begin{equation}\label{eqn:val2}
\mu_q \ge \min(2 \epsilon, \; \delta+2\ord_q(2)) \qquad 
\text{with equality if 
$\; 2 \epsilon \ne \delta+2\ord_q(2)$}.
\end{equation}

\medskip

We deal first with the case where $q=2 \mid d$ (so that $\epsilon=0$). By \eqref{eqn:val1},
we obtain that $(\epsilon,\delta) \equiv (0,1-\nu_2) \mod{\ell}$,
and, by definition, $T_2=\{(0,1-\nu_2)\}$ establishing our claim. 
Next we suppose that $q=2 \nmid d$ (in which case $\nu_2=0$):
\begin{itemize}
\item If $2 \epsilon=\delta+2$ then, from \eqref{eqn:val1} and the fact that $\ell \ge 5$,
we obtain
$(\epsilon,\delta) \equiv (1,0) \mod{\ell}$.
\item If $2 \epsilon> \delta+2$ then, from \eqref{eqn:val2}, we have 
$\mu_2=\delta+2$, so from \eqref{eqn:val1} we obtain
$(\epsilon,\delta) \equiv (3-\mu_2,\mu_2-2) \mod{\ell}$.
\item If $2 \epsilon < \delta+2$ then, from \eqref{eqn:val2}, we
have $\mu_2=2 \epsilon$, so from \eqref{eqn:val1} we obtain
$(\epsilon,\delta) \equiv (\mu_2/2,1-\mu_2/2) \mod{\ell}$.
\end{itemize} 

Next, let us next consider odd $q \mid d$ (whereby we have that  $\mu_q=0$). 
From \eqref{eqn:val2}, it follows that either $\epsilon=0$
or $\delta=0$. From \eqref{eqn:val1}, we obtain
$(\epsilon,\delta) \equiv (0, -\nu_q) $
or $(-\nu_q,0) \mod{\ell}$ as required.

Finally we consider odd $q \mid (d^2-1)$ (so $\nu_q=0$):
\begin{itemize}
\item If $2 \epsilon=\delta$ then, from \eqref{eqn:val1} and the fact that $\ell \ge 5$,
we obtain
$(\epsilon,\delta) \equiv (0,0) \mod{\ell}$.
\item If $2 \epsilon> \delta$ then, from \eqref{eqn:val2}, we have 
$\mu_q=\delta$, so from \eqref{eqn:val1} we obtain
$(\epsilon,\delta) \equiv (-\mu_q,\mu_q) \mod{\ell}$.
\item If $2 \epsilon < \delta$ then, from \eqref{eqn:val2}, we
have $\mu_q=2 \epsilon$, so from \eqref{eqn:val1} we obtain
$(\epsilon,\delta) \equiv (\mu_q/2,-\mu_q/2) \mod{\ell}$.
\end{itemize} 
\end{proof}

From \eqref{eqn:descent} and \eqref{eqn:identity}, we deduce the
following ternary equation
\begin{equation}\label{eqn:ternary}
4\beta y_2^\ell-\alpha^2 y_1^{2\ell}=d^2-1.
\end{equation}
We need to solve this for each possible $(\alpha,\beta) \in \cA_d$
with $2 \le d \le 50$ and $y_1$, $y_2$ integers. Clearing denominators
and dividing by the greatest common divisor of the coefficients we
can rewrite this as
\begin{equation}\label{eqn:ternary2}
r y_2^\ell- s y_1^{2 \ell}=t
\end{equation}
where $r$, $s$, $t$ are positive integers and $\gcd(r,s,t)=1$.

\section{Linear Forms in $2$ logarithms}\label{sec:linear}

The descent step in the previous section transforms \eqref{eqn:main} into a family of ternary
equations \eqref{eqn:ternary}. In this section, we appeal to lower bounds for linear forms in logarithms
to bound the exponent $\ell$ appearing in these equations. 
We will use a special case of Corollary 2 of Laurent \cite{Lau} (with $m=10$ in  the notation of that paper) :
\begin{prop} \label{prop:Laurent}
Let $\alpha_1$ and $\alpha_2$ be positive real, multiplicatively independent algebraic numbers and 
$\log \alpha_1$, $\log \alpha_2$ be any fixed determinations of their logarithms that are real and positive. Write
$D = [ \Q(\alpha_1,\alpha_2) \, : \, \Q ]$ and
\[
b^\prime = \frac{b_1}{D \log A_2} + \frac{b_2}{D \log A_1}
\]
where $b_1$ and $b_2$ are positive integers and $A_1$ and $A_2$ are real numbers $>1$ such that
\[
\log A_i \geq \max \{ h (\alpha_i), \lvert \log \alpha_i\rvert /D, 1/D \}, \qquad i = 1, 2.
\]
Let $\Lambda = b_2 \log \alpha_2 - b_1 \log \alpha_1$.  Then
\[
\log  \lvert \Lambda\rvert \geq -25.2 D^4 \left( \max \{ \log b^\prime + 0.38, 10/D, 1 \} \right)^2 
\log A_1 \log A_2.
\]
\end{prop}
Here, we have defined, as usual, the \emph{absolute logarithmic height} of an algebraic number $\alpha$ by
\[
h(\alpha) = \frac{1}{d} \left( \log |a| + \sum_{i=1}^d \log \max (1,|\alpha^{(i)}|) \right),
\]
where $a$ is the leading coefficient of the minimal polynomial of $\alpha$ and the $\alpha^{(i)}$ are the 
conjugates of $\alpha$ in $\C$.

\bigskip

In this section, we will assume that $3 \le d \le 50$. In the notation of
the previous section,
$(\alpha,\beta)$ will denote an element of  $\cA_d$ while
$(y_1,y_2)$ denotes an integral solution to \eqref{eqn:ternary}. 
By definition of $\cA_d$,
the rationals $\alpha$ and $\beta$ are both positive. It follows from
\eqref{eqn:ternary} that $y_2>0$.
\begin{lem}\label{lem:multind}
Let $\ell > 1000$. 
Suppose $\lvert y_1 \rvert$, $y_2 \ge 2$
and $y_2 \ne y_1^2$.
Let
\begin{equation}\label{eqn:alpha12}
\alpha_1=4 \beta/\alpha^2 \; \; \mbox{ and } \; \;  \alpha_2=y_1^2/y_2 .
\end{equation}
Then $\alpha_1$ and $\alpha_2$ are positive and multiplicatively
independent. Moreover, writing
\begin{equation}\label{eqn:lambda}
\Lambda = \log {\alpha_1} - \ell \log { \alpha_2 }.
\end{equation}
we have
\begin{equation}\label{eqn:lambdabd}
0 < \Lambda < \frac{d^2-1}{\alpha^2 y_1^{2\ell}}.
\end{equation}
\end{lem}
\begin{proof}
By the observation preceding the statement of the lemma, we know that $\alpha_1$ and $\alpha_2$
are positive.
From \eqref{eqn:ternary}, \eqref{eqn:alpha12}, \eqref{eqn:lambda} and \eqref{eqn:lambdabd}, we have
$$
e^\Lambda-1=\frac{4\beta}{\alpha^2} \cdot \frac{y_2^\ell}{y_1^{2\ell}} \, - \, 1=\frac{d^2-1}{\alpha^2 y_1^{2\ell}}
\, > \, 0,
$$
whence $\Lambda>0$. The second part of the lemma thus follows from the inequality
$e^\Lambda-1>\Lambda$.

It remains to show the multiplicative independence of $\alpha_1$ and $\alpha_2$,
so suppose, for a contradiction, that they are
are multiplicatively dependent.
Thus there exist coprime positive integers $u$ and $v$ such that
$\alpha_1^u=\alpha_2^v$. If $\alpha_1=1$ then $\alpha_2=1$
so that $y_2=y_1^2$ contradicting the hypotheses of the lemma.
Thus $\alpha_1 \ne 1$. Defining
\[
g=\gcd\{\ord_p(\alpha_1) \, : \, \text{$p$ prime}\},
\]
as $\alpha_1 \ne 1$, we necessarily have $g \ne 0$. Clearly $v \mid g$. 
However, from \eqref{eqn:lambda},
\[
\Lambda  =(\log{\alpha_1}) \left( 1 - \ell \frac{\log{\alpha_2}}{\log{\alpha_1}}
\right)
= (\log{\alpha_1}) \left( 1- \ell \frac{u}{v} \right) 
= \lvert \log{\alpha_1} \rvert \cdot \left\lvert 1-\ell\frac{u}{v} \right\rvert .
\]
From \eqref{eqn:lambdabd}, we have
\[
0 < \left\lvert 1-\ell\frac{u}{v} \right\rvert < \frac{d^2-1}{\lvert \log {\alpha_1} \rvert \cdot \alpha^2 y_1^{2\ell}}.
\]
Now the non-zero rational $1-\ell u/v$ has denominator dividing $v$ and hence dividing $g$. Thus,
\[
\frac{1}{g} \le \left\lvert 1-\ell\frac{u}{v} \right\rvert  \, .
\]
Since $\lvert y_1 \rvert \ge 2$, it follows that
\[
4^{\ell} \le y_1^{2\ell} < \frac{(d^2-1) g}{\lvert \log {\alpha_1} \rvert \cdot \alpha^2 } \, , 
\]
and so
\[
\ell < \log \left( \frac{(d^2-1) g}{\lvert \log {\alpha_1} \rvert \cdot \alpha^2 } \right)/\log{4} \, .
\]
We wrote a simple \texttt{Magma} script that computes this bound on $\ell$ 
for the values
of $d$ in the range $3 \le d \le 50$ and the possible pairs $(\alpha,\beta) \in \cA_d$
with corresponding $\alpha_1=4\beta/\alpha^2 \ne 1$. We found that the largest possible
value for the right-hand side of the inequality is $19.09\dots$ corresponding to
$d=50$ and $(\alpha,\beta)=(1/62475, 2499)$. As $\ell>1000$, we have a contradiction
by a wide margin.

In fact, we found only one pair $(\alpha,\beta)$ for which $\alpha_1=1$. This
arises when $d=8$ and $(\alpha,\beta)=(1,1/4)$.
\end{proof}

\begin{lem}\label{lem:y1y2}
Let $A_2=\max\{  y_1^2 , \, y_2 \}$. Under the notation and
assumptions of the previous lemma,
\[
1 \; \le \; \frac{\log{A_2} }{\log { y_1^2 }} \;  \le \;  1.03.
\]
\end{lem}
\begin{proof} 
It is sufficient to show that $\log{y_2}/\log{ y_1^2 } \le 1.03$.
From \eqref{eqn:alpha12}, \eqref{eqn:lambda} and \eqref{eqn:lambdabd}, we have 
$$
\log{\alpha_1} - \ell (\log{  y_1^2 }-\log{ y_2 }) \; < \; \frac{d^2-1}{\alpha^2 \cdot 4^\ell}
$$
where we have used the assumption $\lvert y_1 \rvert \ge 2$. It follows that 
\begin{equation*}
\begin{split}
\frac{\log{y_2}}{\log{y_1^2}} & < \;  1 + \frac{1}{\ell \log{y_1^2}} \left(-\log{\alpha_1}+\frac{(d^2-1)}{\alpha^2 \cdot 4^\ell} \right) \\
& \le \; 1 +  \frac{1}{\ell \log{y_1^2}} \left(\lvert \log{\alpha_1} \rvert+\frac{(d^2-1)}{\alpha^2 \cdot 4^\ell} \right) \\
& < \; 1 +  \frac{1}{1000 \log{4}} \left(\lvert \log{\alpha_1} \rvert+\frac{(d^2-1)}{\alpha^2 \cdot 4^{1000}} \right) \,, 
\end{split}
\end{equation*}
using  the assumptions $\ell>1000$ and $\lvert y_1\rvert \ge 2$. We wrote a \texttt{Magma} script that computed
this upper bound for $\log{ y_1^2 }/\log{  y_2 }$ for all $3 \le d \le 50$ and $(\alpha,\beta) \in \cA_d$. 
The largest value of the upper bound we obtained was $1.02257\dots$, again corresponding to
$d=50$ and $(\alpha,\beta)=(1/62475, 2499)$. This completes the proof.
\end{proof}

We continue under the assumptions of Lemma~\ref{lem:multind}, applying
Proposition~\ref{prop:Laurent} to obtain a bound for the exponent $\ell$.
We let 
\[
A_1=\max\{ H(\alpha_1), \,  e \},
\]
where $H(u/v)$, for coprime integers $u$, $v$ (with $v$ non-zero)
is simply $\max\{\lvert u \rvert, \, \lvert v \rvert\}$. Let $A_2$ be
as in Lemma~\ref{lem:y1y2}. We see, thanks to Lemma~\ref{lem:multind},
that the hypotheses of Proposition~\ref{prop:Laurent} are satisfied
for our choices of $\alpha_1$, $\alpha_2$, $A_1$, $A_2$ with $D=1$. We write
\[
b^\prime=\frac{1}{\log{A_2}}+\frac{\ell}{\log{A_1}} > \frac{1000}{\log{A_1}}
\]
as $\ell >1000$. We checked that the smallest possible value for $1000/\log{A_1}$
for $3 \le d \le 50$ and $(\alpha,\beta) \in \cA_d$ is $31.95\cdots$ arising
from the choice $d=50$ and $(\alpha,\beta)=(1/62475, 2499)$. From Proposition~\ref{prop:Laurent},
\[
-\log{\lvert \Lambda \rvert} < 25.2 \log{A_1}\cdot \log{A_2} \cdot (\log{b^\prime})^2 
\le 25.2 \log{A_1} \cdot \log{A_2}\cdot \log^2 \left( \frac{\ell}{\log{A_1}}+\frac{1}{\log{4}}  \right),
\]
where we have used the fact that $A_2 \ge y_1^2 \ge 4$. Combining this with \eqref{eqn:lambdabd}, we have
\[
\ell \log{y_1^2} \; < \;\log\left(\frac{d^2-1}{\alpha^2}\right)+25.2 \log{A_1} \cdot \log{A_2} \cdot \log^2 \left(\frac{\ell}{\log{A_1}}+\frac{1}{\log{4}}  \right).
\]
Next we divide by $\log{y_1^2}$, making use of the fact that $\log{A_2}/\log{y_1^2}<1.03$ and also that
$\lvert y_1 \rvert \ge 2$, to obtain
\[
\ell \; < \; \frac{1}{\log{4}}\log\left(\frac{d^2-1}{\alpha^2}\right)+26 \log{A_1} \cdot \log^2 \left(\frac{\ell}{\log{A_1}}+\frac{1}{\log{4}}  \right).
\]
The only remaining variable in this inequality is $\ell$. It is a 
straightforward exercise in calculus to deduce a bound on $\ell$ for 
any $d$, $\alpha$ and $\beta$. In fact the largest bound on $\ell$
we obtain for $d$ in our range is $\ell< 2,648,167$. We summarize
the results of this section in the following lemma.

\begin{lem}\label{lem:bound}
Let $3 \le d \le 50$ and $(\alpha,\beta) \in \cA_d$. Let $(y_1,y_2)$
be an integral solution to \eqref{eqn:ternary} with $\lvert y_1 \rvert$,
$y_2\ge 2$ and $y_2 \ne y_1^2$. Then $\ell< 3 \times 10^6$.
\end{lem}

\subsection{Proof of Theorem~\ref{thm:main}: bounding $\ell$}
We have dealt with the cases $\ell=2$ and $d=2$ in 
Sections~\ref{sec:l2} and~\ref{sec:d2} respectively,
and so $\ell \ge 3$ and $3 \le d \le 50$. 
We will deal with $\ell=3$ in Section~\ref{sec:l3},
so suppose $\ell \ge 5$. Lemma~\ref{lem:descent}
provides a finite set $\cA_d$ of pairs $(\alpha,\beta)$
such that for every solution $(x,y)$ of \eqref{eqn:main}
there is a pair $(\alpha,\beta) \in \cA_d$ and integers $(y_1,y_2)$
satisfying \eqref{eqn:descent}, \eqref{eqn:ternary}
and \eqref{eqn:ternary2}.
Lemma~\ref{lem:bound}
tells us that $\ell< 3 \times 10^6$ provided the $\lvert y_1 \rvert$,
$y_2>2$
and $y_2 \ne y_1^2$.
It is easy to determine $(y_1,y_2)$ for which these
conditions fail. Indeed, instead of  \eqref{eqn:ternary}
consider the equivalent \eqref{eqn:ternary2} with integral coefficients.
If $y_2=y_1^2$ then \eqref{eqn:ternary2} reduces to $(r-s)y_1^{2\ell}=t$
which allows us to easily determine the corresponding solutions, and similarly for $y_2=1$,
and for $y_1 \in \{ -1, 0, 1 \}$.  We determined all the solutions $(y_1,y_2)$ where the
hypotheses fail for $3 \le d \le 50$ and checked that none of these leads to a solution
to \eqref{eqn:main} with $x\ge 1$ integral (for the purpose
of proving Theorem~\ref{thm:main}, we are only interested in $x \ge 1$).
Thus we may suppose that the hypotheses of Lemma~\ref{lem:bound} hold
and conclude that $\ell<3 \times 10^6$.

\section{A Criterion for the non-existence of solutions}\label{sec:criterion}

In Section~\ref{sec:descent},
we reduced the problem of solving
equation \eqref{eqn:main} (for $3 \le d \le 50$ 
and prime exponents $\ell \ge 5$) 
to the resolution of a number of equations of the form \eqref{eqn:ternary2}. In 
Section~\ref{sec:linear}, we showed that the exponent $\ell$ is necessarily bounded by $3 \times 10^6$.
In this section, we
will provide a criterion for the non-existence of solutions to \eqref{eqn:ternary2}, given $r$, $s$, $t$ and $\ell$.
\begin{lem}\label{lem:criterion}
Let $\ell \ge 3$ be prime. Let $r$, $s$ and $t$ be positive integers satisfying
$\gcd(r,s,t)=1$. Let $q=2k \ell+1$ be a prime that does
not divide $r$. Define
\begin{equation}\label{eqn:mu}
\mu(\ell,q)=\{ \eta^{2\ell} \; : \; \eta \in \F_q \}
=\{0\} \cup \{ \zeta \in \F_q^* \; : \; \zeta^{k}=1\}
\end{equation}
and
$$
B(\ell,q)=\left\{ \zeta \in \mu(\ell,q) \; : \; ((s \zeta+t)/r)^{2k} \in \{0,1\} \right\} \, .
$$
If $B(\ell,q)=\emptyset$,
then equation~\eqref{eqn:ternary2} 
does not have integral solutions.
\end{lem}
\begin{proof}
Suppose $B(\ell,q)=\emptyset$. 
Let $(y_1,y_2)$ be a solution to \eqref{eqn:ternary2}.
Let $\zeta=\overline{y_1}^{2\ell} \in \mu(\ell,q)$. From \eqref{eqn:ternary2}
we have
\[
(s \zeta+t)/r \equiv y_2^\ell \mod{q}.
\]
Thus
\[
((s \zeta+t)/r)^{2k} \, \equiv \,  y_2^{q-1} \,  \equiv \, 
\text{$0$ or $1$} \mod{q}.
\]
This shows that $\zeta \in B(\ell,q)$ giving a contradiction.
\end{proof}

\noindent \textbf{Remark.} We now provide a heuristic explanation
why Lemma~\ref{lem:criterion} should succeed in proving the non-existence of solutions
to \eqref{eqn:ternary2} provided there are no solutions,
particularly if $\ell$ is large. Observe that $\#\mu(\ell,q)=k+1$.
For $\zeta \in \mu(\ell,q)$, the element $((s\zeta+t)/r)^{2k} \in \F_q$ is 
either $0$ or an $\ell$-th root of unity. Thus the ``probability''
that it belongs to $\{0,1\}$ is $2/(\ell+1)$. It follows
that the ``expected
size'' of $B(\ell,q)$ is $2(k+1)/(\ell+1) \approx 2q/\ell^2$. For large
$\ell$ we expect to find a prime $q=2k\ell+1$ such that $2q/\ell^2$
is tiny and so we likewise expect  that $\# B(\ell,q)=0$.

\subsection{Proof of Theorem~\ref{thm:main}: applying the criterion}
\label{sub:criterion}
We wrote a \texttt{Magma} script which, for each $3 \le d \le 50$, 
and each $(\alpha,\beta) \in \cA_d$ (and corresponding triple
of coefficients $(r,s,t)$), and every prime $5 \le \ell <3\times 10^6$,
systematically searches for a prime $q=(2k \ell+1) \nmid r$ 
with $k \le 1000$ such that
$B(\ell,q)=\emptyset$. If it finds such a $q$ then by Lemma~\ref{lem:criterion}
we know that \eqref{eqn:ternary} has no solutions, and thus there are no solutions
to \eqref{eqn:main} that give rise to the pair $(\alpha,\beta)$ via Lemma~\ref{lem:descent}. The entire time for the computation was roughly 3 hours
on a 2500MHz AMD Opteron.
The criterion systematically failed for all exponents 
$5 \le \ell< 3 \times 10^6$ 
whenever $4\beta=d^2-1$ (equivalently the coefficients of \eqref{eqn:ternary2}
satisfy $r=t$). This failure is unsurprising as equations \eqref{eqn:ternary}
and~\eqref{eqn:ternary2} have the obvious solution $(y_1,y_2)=(0,1)$.
In all cases where $4\beta \ne d^2-1$, the criterion succeeded for all
values of $\ell$ except for  a handful of small values.
There were a total of $224$ quintuples $(d,\ell,r,s,t)$ 
with $r\ne t$ for which the criterion fails.
The largest value of $\ell$ in cases $r \ne t$ for which
the criterion fails is $\ell=19$ with $d=27$, $\alpha=1/7$, $\beta=14/27$,
and corresponding $r=2744$, $s=27$, $t=963144$.

\medskip

At this point, to complete the proof of Theorem \ref{thm:main}, we thus require
another method to handle \eqref{eqn:ternary2} when $r=t$,
and also some new techniques to solve this equation when $r \ne t$,
for the remaining small $\ell$. The first question is addressed in Section~\ref{sec:Frey},
and the second in Section~\ref{sec:adhoc}.

\section{Frey-Hellegouarch Curve for the case $r=t$}\label{sec:Frey}

In practice, we have found that Lemma~\ref{lem:criterion} 
will eliminate all elements $(\alpha,\beta) \in \cA_d$
for any given sufficiently large $\ell$
except when $\beta=(d^2-1)/4$ (which is equivalent to
$r=t$). In this case, equation \eqref{eqn:ternary} has the solution
$(y_1,y_2)=(0,1)$ which causes the criterion 
of Lemma~\ref{lem:criterion} fails;
for this situation, we would like to show that $(y_1,y_2)=(0,1)$
is in fact the only solution. In this section, we will thus focus
on \eqref{eqn:ternary} for $\beta=(d^2-1)/4$, and continue
to suppose that $\ell \ge 5$ is prime. 
It follows from the definition of $\cA_d$
that $\alpha=8/d(d^2-1)$, and moreover that 
this pair $(\alpha,\beta)=(8/d(d^2-1),(d^2-1)/4)$
arises exactly when either $\ord_2(d)=0$ or $3$.
We can rewrite \eqref{eqn:ternary} as
\begin{equation}\label{eqn:prepreFrey}
y_2^\ell-\frac{
64 
}{d^2(d^2-1)^3} \cdot
y_1^{2 \ell}
=1.
\end{equation}
We note from \eqref{eqn:descent} that $y_1$ is even if $\ord_2(d)=0$
and $y_1$ is odd if $\ord_2(d)=3$.
By the conclusion of Lemma~\ref{lem:descent}, we know that $y_1$, $y_2$ are integers. It follows from \eqref{eqn:prepreFrey} that
$S \mid y_1$ where 
\[
\begin{cases}
S=\Rad \left(d(d^2-1) \right)  & \text{if $\ord_2(d)=0$,}\\
S=\Rad_2 \left(d(d^2-1)) \right) & \text{if $\ord_2(d)=3$.}
\end{cases}
\]
Let $y_1=S y_3$. Then, from \eqref{eqn:prepreFrey},
\begin{equation}\label{eqn:preFrey}
y_2^\ell- T y_3^{2\ell}=1
\end{equation}
where 
\[
T=\frac{64 S^{2\ell}}{d^2(d^2-1)^3}. 
\]
In addition to the assumption $\ell \ge 5$, let us further suppose that 
\begin{equation}\label{eqn:assumption1}
2 \ell > \ord_q(d^2 (d^2-1)^3)
\end{equation}
for all odd primes $q$. If $\ord_2(d)=0$, we will also assume  that
\begin{equation}\label{eqn:assumption2}
2 \ell \ge 3\ord_2(d^2-1)-1.
\end{equation}
From assumptions~\eqref{eqn:assumption1} and ~\eqref{eqn:assumption2},
it follows that $T$ is an integer and that $\Rad(T)=S$. If $\ord_2(d)=0$, then $2^5 \mid T$.
If, however, $\ord_2(d)=3$, then $\ord_2(T)=0$ and $2 \nmid y_3 \mid y_1$
so that $2 \mid y_2$. We would like to show that all
solutions to \eqref{eqn:prepreFrey} satisfy $y_1=0$, so suppose
$y_1 \ne 0$ (which implies $y_3 \ne 0$). Clearly $y_2 \ne 0$.
We associate our solution $(y_2,y_3)$ to the Frey--Hellegouarch curve
\[
\begin{cases}
E \; : \; Y^2=X(X+1)(X-T y_3^{2 \ell}) & \text{if $\ord_2(d)=0$,}\\
E \; : \; Y^2=X(X+1)(X+y_2^\ell) & \text{ if $\ord_2(d)=3$.} 
\end{cases}
\]
The condition $y_2 y_3 \ne 0$ ensures that the given Weierstrass model
is smooth.  We apply the recipes of 
Kraus \cite{KrausFermat} which build on modularity of 
elliptic curves due to Wiles,
Breuil, Conrad, Diamond and Taylor \cite{Wiles}, \cite{modularity}, 
on Ribet's level lowering
theorem \cite{Ribet}, and on Mazur's theorem \cite{Mazur}. 
The recipes of Kraus are also
reproduced in \cite[Section 14.1]{Siksek}. In the notation
of that reference,
$E \sim_\ell f$ where $f$ is a weight $2$ newform of level 
\[
N=\begin{cases} 
S & \text{if $\ord_2(d)=0$}\\
2 S & \text{if $\ord_2(d)=3$}.
\end{cases}
\] 
If $f$ is irrational (i.e. the Fourier coefficients of $f$ do not all lie in $\mathbb{Q}$) then we can obtain a sharp bound for $\ell$
as we now explain. 
Let $K$ be the number field generated by the coefficients of $f$.
For a prime $q \nmid N$, write $a_q(f) \in \OO_K$ for the $q$-th coefficient
of $f$. Let
\[
H_q=\{ a \in \Z \cap [-2 \sqrt{q},2\sqrt{q}] \; : \; q+1-a \equiv 0 \mod{4}
\}.
\] 
Let
\[
B_q(f)=q \cdot \norm_{K/\Q}((q+1)^2-a_q(f)^2) \cdot \prod_{a \in H_q} \norm_{K/\Q} (a-a_q(f)).
\] 
If $E \sim_\ell f$ then by \cite[Proposition 9.1]{Siksek},
$\ell \mid B_q(f)$. As $f$ is irrational, there is
a positive density of primes $q \nmid N$ such that $a_q(f) \notin \Q$,
and so $B_q(f) \ne 0$. This means that we obtain a bound for $f$,
which is practice is quite small. We can usually improve on this bound
by choosing a set of primes $\cQ=\{q_1,\dotsc,q_n\}$ all
not dividing $N$ and letting
\[
B_{\cQ}(f)=\gcd(B_q(f) \; : \; q \in \cQ).
\]
If $E \sim_\ell f$ then $\ell \mid B_{\cQ}(f)$.

\begin{lem}\label{lem:criterionFrey}
Let $3 \le d \le 50$ with $\ord_2(d)=0$ or $3$.
Suppose $\ell \ge 5$ is a prime that satisfies \eqref{eqn:assumption1}
for all odd primes $q$.
If $\ord_2(d)=0$, suppose $\ell$ also satisfies \eqref{eqn:assumption2}.
Let $N$ be as above. Suppose for each irrational newform of 
weight $2$ and level $N$ there is a set of primes $\cQ$
not dividing $N$ such that $\ell \nmid B_{\cQ}(f)$.
Suppose for every elliptic curve $F$ of conductor $N$ there
is a prime
$q=2k \ell+1$, $q \nmid N$, such that
\begin{enumerate}[(i)]
\item $B(\ell,q)=\{\overline{0}\}$,
where $B(\ell,q)$ is as in the statement of Lemma~\ref{lem:criterion};
\item $\ell \nmid (a_q(F)^2-4)$.
\end{enumerate}
Then 
\begin{itemize}
\item if $\ord_2(d)=3$ then \eqref{eqn:main} 
has no solutions with 
$(\alpha,\beta)=(8/d(d^2-1), (d^2-1)/2)$ in Lemma~\ref{lem:descent};
\item if $\ord_2(d)=0$ then the only solution
to \eqref{eqn:main} with 
$(\alpha,\beta)=(8/d(d^2-1), (d^2-1)/2)$ in Lemma~\ref{lem:descent}
satisfies $x=-(d+1)/2$.
\end{itemize}
\end{lem}
\begin{proof}
The conclusion of the lemma is immediate if $y_1=0$ in \eqref{eqn:descent}.
Let us thus suppose that  $y_1 \ne 0$ and attempt  to deduce a contradiction.
From the above discussion, there is a newform $f$
of level $N$ such that $E \sim_\ell f$, where $E$
is the Frey--Hellegouarch curve. If $f$ is irrational then
$\ell \mid B_{\cQ}(f)$, which contradicts the hypotheses of
the lemma. Thus $f$ is rational and so $f$ corresponds to 
an elliptic curve $F/\Q$ of conductor $N$. Thus $E \sim_\ell F$.

Suppose (i).  By the proof of Lemma~\ref{lem:criterion} we have 
that $q \mid y_1$. Thus $q \mid y_3$.
It follows that $E$ has multiplicative reduction at $q$. 
Thus $(q+1) \equiv \pm a_q(F) \mod{\ell}$.
As $q \equiv 1 \mod{\ell}$ we obtain $4 \equiv a_q(F)^2 \mod{\ell}$.
This contradicts (ii) and completes the proof. 
\end{proof}

\noindent \textbf{Remark.} 
In this section, we are concerned with equation \eqref{eqn:ternary} with $4\beta=d^2-1$,
or equivalently equation \eqref{eqn:ternary2} with $r=t$. These have the solution
$(y_1,y_2)=(0,1)$. It follows from the proof of Lemma~\ref{lem:criterion}
that $\overline{0} \in B(\ell,q)$ (for any suitable $q$) and thus $B(\ell,q) \ne \emptyset$.
However, in this case, the heuristic remark following the proof of Lemma~\ref{lem:criterion}
leads us to expect $B(\ell,q)=\{\overline{0}\}$ for sufficiently large $\ell$ (and suitable
$q$).

\subsection{Proof of Theorem~\ref{thm:main}: the case $r=t$}
\label{sub:r=t}
We wrote a \texttt{Magma} script which, for each $3 \le d \le 50$
with $\ord_2(d)=0$ or $3$,
computes the newforms of weight $2$, level $N$.
Our script take $\cQ$ to be set of primes $<100$ that do not divide $N$,
and computes $B_\cQ(f)$ for each irrational eigenform $f$ at level $N$.
These unsurprisingly are all non-zero. For every prime 
$5 \le \ell< 3 \times 10^6$
that does not divide any of the $B_\cQ(f)$, and satisfies inequality
\eqref{eqn:assumption1}, and also inequality \eqref{eqn:assumption2}
if $\ord_2(d)=0$, and for every isogeny class of elliptic curves $F$
of conductor $N$,
the script systematically searches for a prime $q=(2k \ell+1) \nmid r$ 
with $k \le 1000$ such that conditions (i) and (ii) of 
Lemma~\ref{lem:criterionFrey} hold. If it finds such a $q$ 
we know that there are no solutions
to \eqref{eqn:main} that give rise to the pair $(\alpha,\beta)=(8/d(d^2-1), (d^2-1)/2)$ via Lemma~\ref{lem:criterionFrey}. The entire time for the computation was 
roughly 2.5 hours
on a 2500MHz AMD Opteron.  
In all cases the criterion succeeded for all
values of $\ell$ except for  a handful of small values.
There were a total of $53$ quintuples $(d,\ell,r,s,t)$ 
with $r=t$ for which either $\ell$
does not satisfy the inequalities \eqref{eqn:assumption1},
\eqref{eqn:assumption2}, or it divides $B_\cQ(f)$ for some
irrational eigenform, or for which the script
did not find a suitable $q$ that satisfies (i), (ii).
The largest value of $\ell$ among the $53$
quintuples
is $\ell=19$:  with 
$d=37$, 
$r=t=54762310872$, $s=1$, and with  $d=40$,  
$r=t=102208119975$, $s=1$.

\section{Descent for $\ell =3 $}\label{sec:l3}

In this section we modify the approach of Section~\ref{sec:descent}
to deal with equation \eqref{eqn:main} with exponent $\ell=3$. 

For an integer $m$, we denote by $[m]$ the element in $\{0,1,2\}$  
such that $m \equiv [m] \mod{3}$.
For a prime
$q$ we let $\mu_q$ and $\nu_q$ be as in \eqref{eqn:munu}.
For each prime $q$, we define a finite subset $T_q \subset \{(m,n)  \, : \, 
m,\, n \in \{0,1,2\}  \}$. 
\begin{itemize}
\item If $q \nmid d(d^2-1)$ then let $T_q=\{(0,0)\}$.
\item For $q=2$ we let
\[
T_2=\begin{cases}
\left\{(0,\, [1-\nu_2])\right\} & \text{if $2 \mid d$} \\
\left\{(1,0), \; (0,1),\; (2,2) \right\} & \text{if $2 \nmid d$ and $\mu_2 \ge 4$}. \\
\left\{(1,0), \; (0,1) \right\} & \text{if $2 \nmid d$ and $\mu_2 =3$}. \\
\end{cases}
\]
\item For odd $q \mid d$, let
\[
T_q=\{([-\nu_q],\, 0), \; (0,\, [-\nu_q])\}.
\]
\item For odd $q \mid (d^2-1)$, let 
\[
T_q=\begin{cases}
\{(0,0), \; (1,2), \; (2,1) \} & \text{if $\mu_q \ge 2$}\\
\{(0,0), \; (2,1) \} & \text{if $\mu_q =1$.}\\
\end{cases}
 \]
\end{itemize}

Let $\cA_d$ be the set of pairs of positive integers $(\alpha,\beta)$
such that
$(\ord_q(\alpha),\ord_q(\beta)) \in T_q$ for all primes $q$.

\begin{lem}\label{lem:cubesdescent}
Let $(x,y)$ be a solution to \eqref{eqn:main}
where $\ell =3$ a prime.
Then there are integers $y_1$, $y_2$ and a pair 
$(\alpha,\beta) \in \cA_d$ such that \eqref{eqn:descent} holds.
\end{lem}

\begin{proof}
The proof is an easy adaptation of the proof of Lemma~\ref{lem:descent}. We omit
the details.
\end{proof}

\subsection{Proof of Theorem~\ref{thm:main}: descent for $\ell=3$}
\label{sub:l3}
From this lemma and \eqref{eqn:identity} we reduce the resolution
of \eqref{eqn:main} with $\ell=3$ to solving a number of equations of the 
form \eqref{eqn:ternary}. These can be transformed by 
clearing denominators
and dividing by the greatest common divisor of the coefficients into 
equations of the form \eqref{eqn:ternary2}
where $r$, $s$, $t$ are positive integers and $\gcd(r,s,t)=1$.
An implementation of above procedure leaves us with
$942$ quintuples $(d,\ell,r,s,t)$ with 
$\ell=3$.

\medskip

We emphasize in passing the difference between the approach
of Section~\ref{sec:descent} and that of this section;
the former gives the same set of triples $(r,s,t)$ 
for all exponents $\ell \ge 5$,
whereas the latter gives a possibly different set of triples $(r,s,t)$ 
for $\ell=3$.

\section{Completing the proof of Theorem~\ref{thm:main}}
\label{sec:adhoc}

Looking back at \ref{sub:criterion}, \ref{sub:r=t} and \ref{sub:l3}
we see that, to complete the proof of 
Theorem~\ref{thm:main}, we need to solve $224+53+942=1219$ equations of
the form \eqref{eqn:ternary2} with $r$, $s$ and $t$
positive integers and $\gcd(r,s,t)=1$. In the second
column of Table~\ref{table:end}
we give a breakdown of these equations according to the exponent $\ell$.
In what follows we look at three methods of eliminating or solving
these equations. 

\begin{table}[h]
\centering
{
\begin{tabu}{|c|c|c|c|}
\hline
\multirow{3}{*}{Exponent $\ell$} & original number  & number surviving  & number surviving\\
				 &  of equations \eqref{eqn:ternary2} &   after local & after further\\
				 &  with exponent $\ell$ & solubility tests  & descent\\
\hline\hline
$3$ & $942$ &  $393$ & $223$  \\
\hline
$5$ & $179$ &  $63$ & $3$ \\
\hline
$7$ & $77$  &  $35$ & $0$ \\
\hline
$11$ & $10$ &  $7$ &  $0$\\
\hline
$13$ & $5$  &  $4$ & $0$  \\
\hline
$17$ & $3$  &   $2$ & $0$  \\
\hline
$19$ & $3$  & $3$ &  $0$  \\
\hline\hline
Total & $1219$ & $507$ & $226$ \\
\hline\hline
\end{tabu}
}
\vskip1ex
\caption{In Sections \ref{sub:criterion}, \ref{sub:r=t} and \ref{sub:l3}
we have reduced the proof of 
Theorem~\ref{thm:main} to the resolution of $1219$
equations of the form \eqref{eqn:ternary2}. The first second column
gives a breakdown of this number according to the exponent $\ell$.
The third column gives the number of these equations surviving the 
local solubility tests of Section \ref{sub:locsol}, and the fourth
column gives the number that also survive the further descent
of Section \ref{sub:furtherdesc}.}
\label{table:end}
\end{table}


\subsection{Local Solubility}\label{sub:locsol}
Recall that $\gcd(r,s,t)=1$ in \eqref{eqn:ternary2}.
Write $g=\Rad(\gcd(r,t))$ and suppose that $g > 1$. Then $g \mid y_1$,
and we can write $y_1=g y_1^\prime$, and thus
\[
r y_2^\ell- s g^{2 \ell} {y_1^\prime}^{2\ell}=t. 
\]
Now we may remove a factor of $g$ from the coefficients 
to obtain
\[
r^\prime {y_2}^\ell - s^\prime {y_1^\prime}^{2\ell}=t^\prime,
\]
where $t^\prime=t/g<t$. Likewise, if $h=\gcd(s,t)>1$, we obtain
an equation
\[
r^\prime {y_2^\prime}^\ell - s^\prime {y_1}^{2\ell}=t^\prime,
\]
Likewise 
where $t^\prime=t/h<t$. We apply these operations repeatedly
until we arrive at an equation of the form
\begin{equation}\label{eqn:predescent}
R \rho^\ell-S \sigma^{2\ell}=T
\end{equation}
where $R$, $S$, $T$ are pairwise coprime.
A necessary condition for the existence of solutions is that
for any odd prime $q \mid R$, the residue $-ST$ modulo $q$ is a square.
Besides this simple test we check for local solubility
at the primes dividing $R$, $S$, $T$, and the primes $q \le 19$.
We subjected all of the $1219$ equations to these local
solubility tests. These have allowed us to eliminate $712$
equations, leaving $507$ equations. A breakdown of these
according to the exponent $\ell$ is given in the third column
of Table~\ref{table:end}.

\subsection{A Further Descent}\label{sub:furtherdesc}
If local solubility fails to rule out solutions 
then we carry out a descent to
do so. Specifically, let
\[
S^\prime=\prod_{\text{$\ord_q(S)$ is odd}} q.
\]
Thus $S S^\prime=v^2$. Write $R S^\prime=u$ and $T S^\prime=m n^2$
with $m$ squarefree. We may now rewrite \eqref{eqn:predescent}
as
\[
(v \sigma^\ell+n \sqrt{-m})(v \sigma^\ell-n \sqrt{-m})=u \rho^\ell. 
\]
Let $K=\Q(\sqrt{-m})$ and $\OO$ be its ring of integers.
Let $\sS$ be the prime ideals of $\OO$ that
divide $u$ or $2n \sqrt{-m}$. Clearly 
$(v\sigma^\ell+n \sqrt{-m}) {K^*}^\ell$ belongs to the ``$\ell$-Selmer group''
\[
K(\sS,\ell)=\{\epsilon \in K^*/{K^*}^\ell \; : \; 
\text{$\ord_\mP(\epsilon) \equiv 0 \mod{\ell}$ for all $\mP \notin \sS$}
\}.
\]
This is an $\F_\ell$-vector space of finite dimension and,
for a given $\ell$, easy
to compute from class group and unit group information (see
\cite[Proof of Proposition VIII.1.6]{Silverman}).
Let
\[
\cE=\{ \epsilon \in K(\sS,\ell) \; : \; \Norm(\epsilon)/u \in {\Q^*}^\ell \}.
\]
It follows that
\begin{equation}\label{eqn:furtherdescent}
v \sigma^\ell+n \sqrt{-m}=\epsilon \eta^\ell,
\end{equation}
where $\eta \in K^*$ and $\epsilon \in \cE$.

\begin{lem}\label{lem:valuative}
Let $\fq$ be a prime ideal of $K$. Suppose one of the following
holds:
\begin{enumerate}
\item[(i)] $\ord_\fq(v)$, $\ord_\fq(n\sqrt{-m})$, $\ord_\fq(\epsilon)$
are pairwise distinct modulo $\ell$;
\item[(ii)] $\ord_\fq(2v)$, 
$\ord_\fq(\epsilon)$, $\ord_\fq(\overline{\epsilon})$
are pairwise distinct modulo $\ell$;
\item[(iii)] $\ord_\fq(2 n \sqrt{-m})$, 
$\ord_\fq(\epsilon)$, $\ord_\fq(\overline{\epsilon})$
are pairwise distinct modulo $\ell$.
\end{enumerate}
Then there is no $\sigma \in \Z$
and $\eta \in K$ satisfying \eqref{eqn:furtherdescent}.
\end{lem}
\begin{proof}
Suppose (i) holds. Then the three terms in \eqref{eqn:furtherdescent}
have pairwise distinct valuations, so \eqref{eqn:furtherdescent}
is impossible $\fq$-adically.
If (ii) or (iii), then we apply the same idea to
\[
2v \sigma^\ell= \epsilon\, \eta^\ell+\overline{\epsilon}\, \overline{\eta}^\ell,
\qquad
2 n \sqrt{-m}= \epsilon\, \eta^\ell-\overline{\epsilon}\, \overline{\eta}^\ell \, ,
\]
which follow from \eqref{eqn:furtherdescent}, and its conjugate
equation.
\end{proof}

\begin{lem}\label{lem:furtherdescent}
Let $q=2k \ell+1$ be a prime. 
Suppose $q\OO=\fq_1 \fq_2$ where $\fq_1$, $\fq_2$
are distinct, and such that $\ord_{\fq_j}(\epsilon)=0$
for $j=1$, $2$. Let 
\[
\chi(\ell,q)=\{ \eta^\ell \; : \; \eta \in \F_q \}.
\]
Let
\[
C(\ell,q)=\{\zeta \in \chi(\ell,q) \; : \;
((v \zeta+n\sqrt{-m})/\epsilon)^{2k} \equiv \text{$0$ or $1 \mod{\fq_j}$
for $j=1$, $2$}\}.
\]
Suppose $C(\ell,q)=\emptyset$. 
Then there is no $\sigma \in \Z$
and $\eta \in K$ satisfying \eqref{eqn:furtherdescent}.
\end{lem}
\begin{proof}
The proof is a straightforward modification
of the proof of Lemma~\ref{lem:criterion}.
\end{proof}
We have found Lemmata~\ref{lem:valuative} and~\ref{lem:furtherdescent}
useful in eliminating many, and often all, $\epsilon \in \cE$.
Of course if they succeed in eliminating all $\epsilon \in \cE$
then we know that \eqref{eqn:predescent} has no solutions,
and so the same would be true for \eqref{eqn:ternary2}.
Of course, when $r=t$, equation \eqref{eqn:ternary2}
always has a solution, namely $(y_1,y_2)=(0,1)$. For
$r=t$, the reduction process in \ref{sub:locsol}
 leads to equation \eqref{eqn:predescent} with $R=T=1$.
The solution $(y_1,y_2)=(0,1)$ to \eqref{eqn:ternary2}
corresponds to the solution
$(\rho,\sigma)=(1,0)$ in \eqref{eqn:predescent}. It follows
from \eqref{eqn:furtherdescent} that 
$n\sqrt{-m} {K^*}^\ell \in \cE$. Naturally,
Lemma~\ref{lem:valuative} and Lemma~\ref{lem:furtherdescent}
do not eliminate the case $\epsilon=n\sqrt{-m}$ since 
equation \eqref{eqn:furtherdescent} has the solution
with $\sigma=0$ and $\eta=1$. In this case, our 
interest is in showing that this is the only solution.
\begin{lem}\label{lem:furtherdescent2}
Suppose 
\begin{enumerate}
\item[(i)] $\ord_\fq(n\sqrt{-m})<\ell$ for all prime ideals $\fq$ of $\OO$;
\item[(ii)] the polynomial $X^\ell+(d-X)^\ell-2$ has no roots
in $\OO$ for $d=1$, $-1$, $-2$;
\item[(iii)] the only root of the polynomial $X^\ell+(2-X)^\ell-2$ 
in $\OO$ is $X=1$.
\end{enumerate}
Then, for $\epsilon=n\sqrt{-m}$,
 the only solution to \eqref{eqn:furtherdescent}
with $\sigma \in \Z$ and $\eta \in K$ is $\sigma=0$
and $\eta=1$.
\end{lem}
\begin{proof}
Let $\epsilon=n\sqrt{-m}$ and suppose $\sigma \in \Z$
and $\eta \in K$ is a solution to \eqref{eqn:furtherdescent}.
Note that the left-hand side of \eqref{eqn:furtherdescent}
belongs to $\OO$, and from (i), we deduce that $\eta \in \OO$.
Now substracting \eqref{eqn:furtherdescent} from its conjugate
and dividing by $n \sqrt{-m}$ leads to the equation
\[
\eta^\ell+\overline{\eta}^\ell=2.
\]
We deduce that the rational integer $\eta+\overline{\eta}$ divides $2$
and hence $\eta+\overline{\eta}=d$ where $d=\pm 1$, $\pm 2$.
Thus $\eta$ is a root of $X^\ell+(d-X)^\ell-2$ for one of 
these values of $d$. By (ii), (iii) it follows that $d=2$
and $\eta=1$. From \eqref{eqn:furtherdescent} we see that 
$\sigma=0$.
\end{proof}
For each of the $507$ equations \eqref{eqn:ternary2} that
survive the local solubility tests in Section \ref{sub:locsol},
we  computed the set $\cE$ and applied the criteria
in Lemma~\ref{lem:valuative} and
Lemma~\ref{lem:furtherdescent} (the latter with $k \le 1000$)
to eliminate as many of the $\epsilon \in \cE$
as possible.
If the two lemmata succeed in eliminating
all possible values of $\epsilon$ then \eqref{eqn:predescent}
has no solutions, and therefore equation \eqref{eqn:ternary2}
does not have solutions either. If
they succeeded in eliminating all but one value
$\epsilon \in \cE$, and that value is $n \sqrt{-m}$,
then we checked the conditions of Lemma~\ref{lem:furtherdescent2}
which if satisfied allow us to conclude that $\sigma=0$
and therefore $y_1=0$. 
Recall that Theorem~\ref{thm:main} is concerned
with \eqref{eqn:main} with $x \ge 1$.
If $y_1=0$ then $x=-(d+1)/2$
(via \eqref{eqn:descent}) and so we can eliminate
$(r,s,t)$ if 
Lemmata~\ref{lem:valuative},~\ref{lem:furtherdescent} and ~\ref{lem:furtherdescent2}
allow us to conclude that $\sigma=0$.
Using this method, we 
managed to eliminate $281$ of the $507$ equations \eqref{eqn:ternary2},
leaving just $226$ equations. In Table~\ref{table:end}
we provide a breakdown of these according to the the exponent $\ell$.

\subsection{A Thue Approach}
Finally, writing $\tau=\sigma^2$ in \eqref{eqn:predescent} we obtain the (binomial)
Thue equation
\[
R \rho^\ell-S \tau^\ell=T. 
\]
We solved the remaining $226$ equations using the 
the Thue equation
solver in \texttt{Magma}.
The theory behind this Thue
equation solver is discussed in  \cite[Chapter VII]{Smart}. 
As we see from Table~\ref{table:end}, we are left with the problem  of
solving $223$ Thue equations of degree $3$, and three
Thue equations of degree $5$.
Working backwards from these solutions, we obtained
precisely six solutions to \eqref{eqn:main}
with $x \ge 1$. These are
\begin{gather*}
3^3+4^3+5^3=6^3, \qquad 11^3+12^3+13^3+14^3=20^3,\\
3^3+4^3+5^3+\cdots + 22^3=40^3,\\
15^3+16^3+17^3+\cdots+34^3=70^3,\\
6^3+7^3+8^3+\cdots+30^3=60^3,\\
291^3+292^3+293^3+\cdots+339^3=1115^3.
\end{gather*}
Noting that these solutions are in Table~\ref{table:solutions},
this completes the proof of Theorem~\ref{thm:main}.


\end{document}